\documentclass[11pt]{article}
\usepackage{graphics}
\usepackage[demo]{graphicx}
\usepackage{epsfig}
\usepackage{pstricks}
\usepackage[normal]{subfigure}
\usepackage[latin1]{inputenc}
\usepackage[english,francais]{babel}
\usepackage{relsize,exscale}
\usepackage{makeidx}
\usepackage{enumitem}
\usepackage{amsfonts,amssymb,amsmath}
\usepackage{graphicx}
\usepackage{color}
\usepackage{multirow}
\usepackage{mathrsfs}
\usepackage[normalem]{ulem}
\usepackage{cancel}
\newenvironment{prooff}{{\it Proof :}}{\hfill\rule{2mm}{2mm}\vskip3mm\par}
\newtheorem{theorem}{Theorem}[section]
\newtheorem{lemma}[theorem]{Lemma}

\newtheorem{corollary}[theorem]{Corollary}
\newtheorem{e-definition}[theorem]{Definition\rm}

\usepackage{color}
\definecolor{dred}{rgb}{0.92,0,0}
\definecolor{dgreen}{rgb}{0,0.92,0}
\definecolor{dblue}{rgb}{0,0,0.92}
\definecolor{dyellow}{rgb}{0.95,0.95,0}

\newcommand{\R}{\mathbb{R}}
\newcommand{\N}{\mathbb{N}}
\def\D{\displaystyle}
\newcommand{\hs}{\hspace{0.1cm}}

\newcommand{\sa}{\\ [0.2cm]}
\usepackage{multirow}
\graphicspath{
{./Figures/}
{./}
}
\title{A new probabilistic interpretation of Bramble-Hilbert lemma}
\author{Jo\"el Chaskalovic \thanks{D'Alembert,
Sorbonne University, Paris, France, (\emph{email}: jch1826@gmail.com)}
\qquad
Franck Assous
\thanks{
Ariel University, 40700 Ariel, Isra\"el, (\emph{email}: franckassous55@gmail.com).}
}
\date{}
\begin{document}
\maketitle
\selectlanguage{english}
\begin{abstract}
\noindent The aim of this paper is to provide new perspectives on relative finite element accuracy which is usually based on the asymptotic speed of convergence comparison when the mesh size $h$ goes to zero. Starting from a geometrical reading of the error estimate due to Bramble-Hilbert lemma, we derive two probability distributions that estimate the relative accuracy, considered as a random variable, between two Lagrange finite elements $P_k$ and $P_m$, ($k < m$). We establish mathematical properties of these probabilistic distributions and we get new insights which, among others, show that $P_k$ or $P_m$ is more likely accurate than the other, depending on the value of the mesh size $h$.
\end{abstract}
\noindent {\em keywords}: Error estimates, Finite elements, Bramble-Hilbert lemma, Probability.
%
%
\section{Introduction}\label{intro}
\noindent The past decades have seen the development of finite element error estimates due to their influence on improving both accuracy and reliability in scientific computing. \sa
However, in these error estimates, an unknown constant is involved which depends, among others, on the basis functions of the considered finite element and on a given semi-norm of the exact solution one wants to approximate. Moreover, error estimates are only upper bounds of the approximation error yielding that the precise value of the approximation error is generally unknown. \sa
Moreover, due to quantitative uncertainties which are generated in the process of the mesh generator and, as a consequence, in the corresponding approximation too, it gave us the idea of considering the approximation error as a random variable. \sa
Therefore, we were able to evaluate the probability of the difference between two approximation errors corresponding to two different finite elements, and then, we got a probabilistic way to compare the relative accuracy between these two finite elements.\sa
The paper is organized as follows. We recall in Section \ref{Geo_and_Proba} the mathematical problem we consider and a corollary of Bramble-Hilbert lemma to propose a geometrical interpretation of the error estimate which appears in this lemma. In Section \ref{models} we derive two probability distributions to interpret and estimate the relative accuracy, considered as a random variable, between two Lagrange finite elements $P_k$ and $P_m, (k<m)$.  Several mathematical properties of these probabilistic distributions are established in Section \ref{probabilistic_law_properties}. Concluding remarks follow.
\section{The problem model and a geometrical interpretation of an error estimate}\label{Geo_and_Proba}
\begin{figure}[h]
  \centering
  \includegraphics[width=10cm]{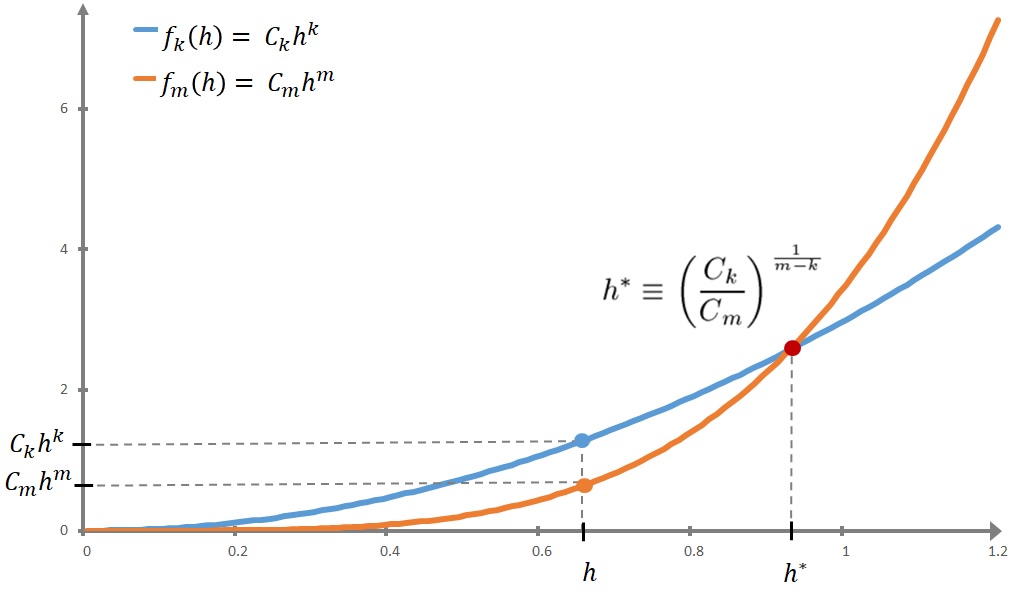}
  \caption{Curves $f_k$ and $f_m$ and existence domain of $\|u^{(i)}_h-u\|_{1,\Omega}, i=k \mbox{ or } i=m$.} \label{Puissance}
\end{figure}

\noindent Let $\Omega$ be an open bounded, and non empty subset of $\R^{n}$ and $\Gamma$ its boundary which we assumed to be $C^1- $piecewise, and let $u$ be the solution to the second order elliptic variational formulation:
\begin{equation}\label{VP}
\textbf{(VP}\textbf{)} \hspace{0.2cm} \left\{
\begin{array}{l}
\mbox{Find } u \in   V \mbox{ solution to:} \\[0.1cm]
a(u,v) = l(v), \quad\forall v \in V,
\end{array}
\right.
\end{equation}
where $V$ is a given Hilbert space endowed with a norm $ \left\|.\right\|_{V}$, $a(\cdot,\cdot)$ is a bilinear, continuous and $V-$elliptic form defined on $V \times V$, and $l(\cdot)$ a linear continuous form defined on~$V$.\sa
Classically, variational problem \textbf{(VP)} has one and only solution $u \in V$ (see for example \cite{ChaskaPDE}). In this paper and for simplicity, we will restrict ourselves to the case where $V$ is a usual Sobolev space of distributions. \sa
Let us also consider an approximation $u_{h}$ of $u$, solution  to the approximate variational formulation:
\begin{equation}\label{VP_h}
\textbf{(VP}\textbf{)}_{h} \hspace{0.2cm} \left\{
\begin{array}{l}
\mbox{Find } u_{h} \in   V_h \mbox{ solution to:} \\[0.1cm]
a(u_{h},v_{h}) = l(v_{h}),\quad \forall v_{h} \in V_h, 
\end{array}
\right.
\end{equation}
where $V_h$ is a given finite-dimensional subset of $V$. \sa
To state a corollary of Bramble-Hilbert's lemma and a corresponding error estimate, we follow \cite{RaTho82} or \cite{Ciarlet}, and we assume that $\Omega$ is exactly recovered by a mesh ${\mathcal T}_h$ composed by $N_K$ n-simplexes $K_{\mu}, (1 \leq \mu \leq N_K),$ which respect classical rules of regular discretization, (see for example \cite{ChaskaPDE} for the bidimensional case and \cite{RaTho82} in $\R^n$). Moreover, we denote by $P_k(K_{\mu})$ the space of polynomial functions defined on a given n-simplex $K_{\mu}$ of degree less than or equal to $k$, ($k \geq$ 1). \sa
Then, we have the following result: \vspace{0.1 cm}
\begin{lemma}\label{Thm_error_estimate}
Suppose that there exists an integer $k \geq 1$ such that the approximation $u_h$ of $V_h$ is a continuous piecewise function composed by polynomials which belong to $P_k(K_{\mu}), (1\leq \mu\leq  N_K)$. \sa
Then, $u_h$ converges to $u$ in $H^1(\Omega)$:
\begin{equation}
\D\lim_{h\rightarrow 0}\|u_h-u\|_{1,\Omega}=0.
\end{equation}
Moreover, if the exact solution $u$ belongs to $H^{k+1}(\Omega)$, we have the following error estimate:
\begin{equation}\label{estimation_error}
\|u_h-u\|_{1,\Omega} \hs \leq \hs \mathscr{C}_k\,h^k \, |u|_{k+1,\Omega}\,,
\end{equation}
where $\mathscr{C}_k$ is a positive constant independent of $h$, $\|.\|_{1,\Omega}$ the classical norm in $H^1(\Omega)$ and $|.|_{k+1,\Omega}$ denotes the semi-norm in $H^{k+1}(\Omega)$.
\end{lemma}
%
Let us now consider two families of Lagrange finite elements $P_k$ and $P_m$ corresponding to a set of values $(k,m)\in \N^2$ such that $0 < k < m$. \\[0.1cm]
The two corresponding inequalities given by (\ref{estimation_error}), assuming that the solution $u$ to \textbf{(VP)} belongs to $H^{m+1}(\Omega)$, are:
\vspace{-0.2cm}
\begin{eqnarray}
\|u^{(k)}_h-u\|_{1,\Omega} \hs & \leq & \hs \mathscr{C}_k h^{k}\, |u|_{k+1,\Omega}, \label{Constante_01} \\
\|u^{(m)}_h\hspace{-0.09cm}-u\|_{1,\Omega} \hs & \leq & \hs \mathscr{C}_m h^{m}\, |u|_{m+1,\Omega}\,, \label{Constante_02}
\end{eqnarray}
where $u^{(k)}_h$ and $u^{(m)}_h$ respectively denotes the $P_k$ and $P_m$ Lagrange finite element approximations of $u$.\\[0.2cm]
Now, if one considers a given mesh for the finite element of $P_m$ which would contains whose of $P_k$ then, for the particular class of problems where $\textbf{(VP)}$ is equivalent to a minimization formulation $\textbf{(MP)}$ (see for example \cite{ChaskaPDE}), one can show that the approximation error of $P_m$ is always lower than those of $P_k$, and $P_m$ is more accurate than $P_m$ for all values of the mesh size $h$ corresponding to the largest diameter in the mesh ${\mathcal T}_h$.\sa
Then, for a given mesh size value of $h$, we consider two independent meshes for $P_k$ and $P_m$ built be a mesh generator. So, usually, to compare the relative accuracy between these two finite elements, one asymptotically considers inequalities (\ref{Constante_01}) and (\ref{Constante_02}) to conclude that, when $h$ goes to zero, $P_m$ finite element is more accurate that $P_k$, as $h^m$ goes faster to zero than $h^k$. \sa
However, for any application $h$ has a static fixed value and this way of comparison is not valid anymore. Therefore, our point of view will be to determine the relative accuracy between two finite elements $P_k$ and $P_m, (k<m)$, for any given value of $h$ for which two independent meshes have to be considered.\sa
To this end, let us set:
\begin{equation}
C_k = \mathscr{C}_k |u|_{k+1,\Omega} \mbox{ and } C_m = \mathscr{C}_m |u|_{m+1,\Omega}.
\end{equation}
Therefore, instead of (\ref{Constante_01}) and (\ref{Constante_02}), we consider in the sequel the two next inequalities:
\begin{eqnarray}
\|u^{(k)}_h-u\|_{1,\Omega} \hs & \leq & \hs C_k h^{k}, \label{Constante_01_2} \\
\|u^{(m)}_h\hspace{-0.09cm}-u\|_{1,\Omega} \hs & \leq & \hs C_m h^{m}. \label{Constante_02_2}
\end{eqnarray}
Then, let us remark that inequalities (\ref{Constante_01_2}) and (\ref{Constante_02_2}) show that the two polynomial curves defined by $f_k(h)\equiv C_k h^k$ and $f_m(h)\equiv C_m h^m$ play a critical role regarding the values of the two norms $\|u^{(k)}_h-u\|_{1,\Omega}$ and $\|u^{(m)}_h-u\|_{1,\Omega}$.\sa
More precisely, these inequalities indicate that the norm $\|u^{(k)}_h-u\|_{1,\Omega}$, (respectively the norm $\|u^{(m)}_h-u\|_{1,\Omega}$), is below the curve $f_k(h)$, (respectively below the curve $f_m(h)$), (see Figure \ref{Puissance}).\sa
As we are interested in comparing the relative positions of these curves, we introduce their intersection point $h^*$ defined by:
\begin{equation}\label{h*}
\D h^* \equiv\left( \frac{C_k}{C_m}\right)^{\frac{1}{m-k}}.
\end{equation}
Now, as often in numerical analysis, there is no {\em a priori} information to surely or better specify the relative distance between $\|u^{(k)}_h-u\|_{1,\Omega}$, (respectively $\|u^{(m)}_h-u\|_{1,\Omega}$), and the curve $f_k$ or its precise value in the interval $[0, C_k h^k]$,  (respectively the curve $f_m$ and the interval $[0, C_m h^m]$). \sa
Moreover, we have to deal with finite element methods that return quantitative uncertainties in their calculations. This mainly comes from the way the mesh grid generator will process the mesh to compute the approximation $u^{(k)}_h$, leading to a partial non control of the mesh, even for a given maximum mesh size. As a consequence, the corresponding grid is \emph{a priori} random, and the corresponding approximation $u^{(k)}_h$ too. \sa
For all of these reasons, we motivate that a probabilistic approach can provide a coherent framework for modeling quantitative uncertainties in finite element approximations. \sa
This is the purpose of the following section where we will establish two probability distributions which will allowed us to estimate the relative accuracy between two Lagrange finite elements.
\vspace{-0.2cm}
\section{The two probabilistic models for relative finite elements accuracy}\label{models}
\vspace{-0.2cm}
\noindent In this section, we will introduce a convenient probabilistic framework to consider the possible values of the norm $\|u^{(k)}_h-u\|_{1,\Omega}$ as a random variable defined as follows:
\begin{itemize}
\item A {\em random trial} corresponds to the grid constitution and the associated approximation $u^{(k)}_h$.
\item The probability space ${\bf\Omega}$ contains therefore all the possible results for a given random trial, namely, all of the possible grids that the mesh generator may processed, or equivalently, all of the corresponding associated approximations $u^{(k)}_h$.
\end{itemize}
Then, for a fixed value of $k$, we define by $X^{(k)}$ the random variable as follows:
\begin{eqnarray}
X^{(k)} : & {\bf\Omega} & \hspace{0.1cm}\rightarrow \hspace{0.2cm}[0,C_k h^k] \noindent \\
& \boldsymbol{\omega}\equiv u^{(k)}_h & \hspace{0.1cm} \mapsto \hspace{0.2cm}\D X^{(k)}(\boldsymbol{\omega}) = X^{(k)}(u^{(k)}_h) = \|u^{(k)}_h-u\|_{1,\Omega}. \label{Def_Xi_h}
\end{eqnarray}
In the sequel, for simplicity, we will set: $X^{(k)}(u^{(k)}_h)\equiv X^{(k)}(h)$. \sa
Now, regarding the absence of information concerning the more likely or less likely values of the norm $\|u^{(k)}_h-u\|_{1,\Omega}$ in the interval $[0, C_k h^k]$, we will assume that the random variable $X^{(k)}(h)$ has a uniform distribution on the interval $[0, C_k h^k]$.\sa
So, our interest is to evaluate the probability of the event
\begin{equation}\label{objectif}
\left\{\|u^{(m)}_h-u\|_{1,\Omega} \leq \|u^{(k)}_h-u\|_{1,\Omega}\right\} \equiv \left\{X^{(m)}(h) \leq X^{(k)}(h)\right\},
\end{equation}
which will allow us to estimate the relative accuracy between two finite elements of order $k$ and $m$, $(k<m)$.\sa
To proceed it, let us now introduce the two random events $A$ and $B$ as follows:
\vspace{-0.2cm}
\begin{eqnarray}
A & \equiv & \left\{\|u^{(m)}_h-u\|_{1,\Omega} \leq \|u^{(k)}_h-u\|_{1,\Omega} \right\}, \label{A}\\
B & \equiv & \left\{\|u^{(k)}_h-u\|_{1,\Omega}\in [C_m h^m,C_k h^k]\right\}. \label{B}
\end{eqnarray}
Then, we have the following lemma:
\begin{lemma}\label{Prob_General}
Let $A$ and $B$ be the events defined by (\ref{A}) and (\ref{B}). Then, we have:
\vspace{-0.2cm}
\begin{equation}\label{Eq_1_5.0}
\D \forall h<h^*: Prob\left\{A\right\} = \frac{Prob\left\{B\right\}}{Prob\left\{B / A\right\}}.
\end{equation}
\end{lemma}
\begin{prooff}
Let us use the following splitting:
\begin{equation}\label{Eq_1}
Prob\left\{A\right\} = Prob\left\{A \cap B\right\} + Prob\left\{A \cap \bar{B}\right\},
\end{equation}
where $\bar{B}$ denotes the opposite event of $B$. \sa
Now, by the definition of the conditional probability we have:
\begin{equation}\label{Eq_1_1}
Prob\left\{A \cap B\right\} = Prob\left\{A / B\right\}.Prob\left\{B\right\} = Prob\left\{B\right\},
\end{equation}
since the probabilistic interpretation of Bramble-Hilbert lemma in the case $h~<~h^*$ corresponds to:
\begin{equation}\label{Eq_1_2}
Prob\left\{A / B\right\} = 1.
\end{equation}
Then, equation (\ref{Eq_1}) can be written as:
\begin{equation}\label{Eq_1_3}
Prob\left\{A\right\} = Prob\left\{B\right\} + Prob\left\{A \cap \bar{B}\right\},
\end{equation}
which can be transformed by the help of the conditional probability as follows:
\begin{equation}\label{Eq_1_4}
Prob\left\{A\right\} = Prob\left\{B\right\} + Prob\left\{\bar{B} / A\right\}.Prob\left\{A\right\},
\end{equation}
or equivalently,
\vspace{-0.2cm}
\begin{equation}\label{Eq_1_5}
\D Prob\left\{A\right\} = \frac{Prob\left\{B\right\}}{1 - Prob\left\{\bar{B} / A\right\}} = \frac{Prob\left\{B\right\}}{Prob\left\{B / A\right\}},
\end{equation}
which corresponds to (\ref{Eq_1_5.0}).
\end{prooff}
Then, we have two options regarding the nature of the dependency between the events $A$ and $B$ which will lead us to get two different distribution laws of probabilities of the event $\D\left\{X^{(m)}(h) \leq X^{(k)}(h)\right\}$. \sa
The next two subsections are devoted to the dependency modeling between $A$ and $B$.
\vspace{-0.2cm}
\subsection{The two steps model}\label{two_steps}
\noindent The first case we will consider states that since, {\em a priori}, no information is available in numerical analysis to consider any kind of dependency between the events $A$ and $B$, we assume in this subsection that these events are independent.
\begin{corollary}
Let $A$ and $B$ be the two events defined by (\ref{A}) and (\ref{B}) and let us assume they are independent. Then, the probability distribution of the event $\D\left\{X^{(m)}(h) \leq X^{(k)}(h)\right\}$ is given by:
\vspace{-0.2cm}
\begin{equation}\label{Heaviside_Prob}
\D Prob\left\{ X^{(m)}(h) \leq X^{(k)}(h)\right\} = \left |
\begin{array}{ll}
\hs 1 & \mbox{ if } \hs 0 < h < h^*, \medskip \\
\hs 0 & \mbox{ if } \hs h> h^*.
\end{array}
\right.
\end{equation}
\end{corollary}
\begin{prooff}
As the events $A$ and $B$ are supposed independent, we have:
\begin{equation}
Prob\left\{A / B\right\} = Prob\left\{A\right\}.
\end{equation}
As a consequence, by lemma \ref{Prob_General} equation (\ref{Eq_1_5}) gives after simplification:
\begin{equation}\label{Eq_1_5.1}
\D Prob\left\{A\right\} = 1, \forall h<h^*.
\end{equation}
With the same kind of arguments, when $h>h^*$ we get:
\begin{equation}\label{Eq_1_5.1.1}
\D Prob\left\{A\right\} = 0, \forall h>h^*.\vspace{-0.3cm}
\end{equation}
\vspace{-0.2cm}
\end{prooff}
Let us now examine the main properties of probabilistic distribution (\ref{Heaviside_Prob}):
\begin{itemize}
\item For any $h$ smaller than $h^*$, $P_m$ finite element is not only asymptotically better than $P_k$ finite element as $h$ becomes small, but they are \emph{almost surely} more accurate for all these values of $h$ such that $h<h^*$.\vspace{0.1cm}
\item For any $h$ greater than $h^*$, $P_k$ finite element becomes \emph{almost surely} more accurate than $P_m$ finite element, even if $k<m$.
\end{itemize}
This last feature upsets the widespread idea regarding the relative accuracy between $P_k$ and $P_m, (k<~\!\!m)$, finite elements. It clearly indicates that there exist cases where $P_m$ finite elements \emph{surely} must be overqualified and a significant reduction of implementation and execution cost can be obtained without a loss of accuracy.\sa
Furthermore, one may expect to get a probabilistic distribution where more variations would appear, as it is in this two steps model, between the probability of the event $\D\left\{ X^{(m)}(h) \leq X^{(k)}(h)\right\}$ and the mesh size $h$. It is certainly due to the assumption we considered regarding the independency between the events $A$ and $B$. \sa
The purpose of the next subsection we will be devoted to relax this assumption by directly computing the probability of the event $\D\left\{ X^{(m)}(h) \leq X^{(k)}(h)\right\}$.
\vspace{-0.2cm}
\subsection{The "sigmoid" model}
\begin{figure}[h]
  \centering
  \includegraphics[width=10cm]{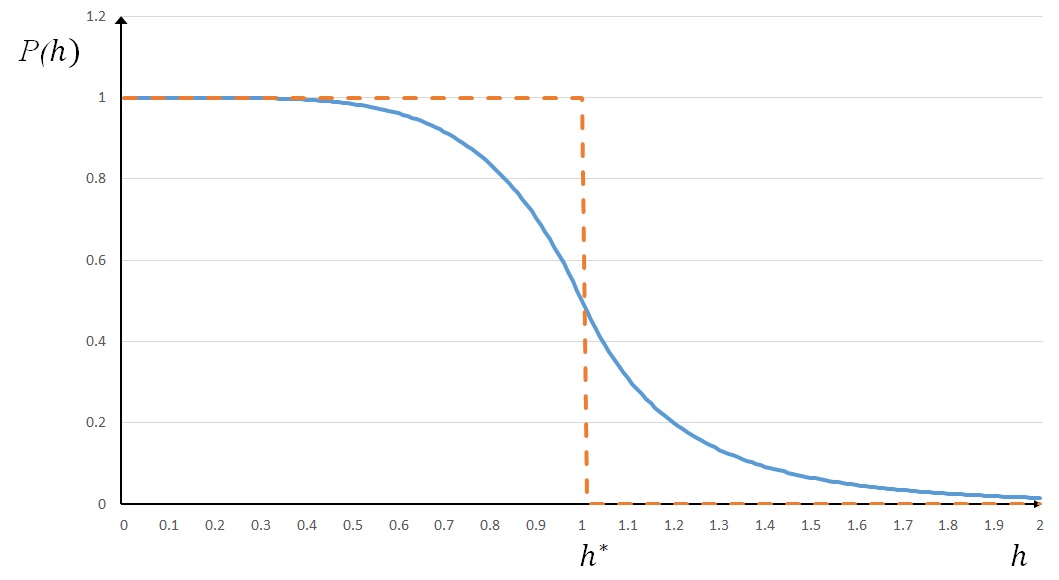}
  \caption{Case $m-k\neq 1$: shape of the sigmoid distribution (\ref{Nonlinear_Prob}) and the two steps corresponding one (\ref{Heaviside_Prob}), $(P(h)\equiv Prob\{X^{(m)}(h) \leq X^{(k)}(h)\})$.} \label{Sigmoid}
\end{figure}

\noindent To avoid the hypothesis of independency between the events $A$ and $B$ defined by (\ref{A}) and (\ref{B}), we will directly evaluate the probability of the event $A$ without considering anymore the splitting we wrote in formula (\ref{Eq_1_3}). \sa
However, we will assume that the two random variables $X^{(i)}(h), (i=k\,\,\mbox{or}\,i=m),$ defined by (\ref{Def_Xi_h}) are independent and uniformly distributed on $[0, C_i h^i], (i=k \mbox{ or } i=m)$.\sa
This is the aim of the following theorem.
\begin{theorem}\label{The_nonlinear_law}
Let $u$ be the solution to the second order variational elliptic problem $\textbf{(VP}\textbf{)}$ defined in (\ref{VP}) and $u^{(i)}_h, (i=k \mbox{ or } i=m, k<m)$, the two corresponding Lagrange finite element $P_i$ approximations, solution to the approximated formulation $\textbf{(VP}\textbf{)}_{h}$ defined by (\ref{VP_h}).\\[0.1cm]
We assume the two corresponding random variables $X^{(i)}(h), (i=k \mbox{ or } i=m),$ defined by (\ref{Def_Xi_h})  are independent and uniformly distributed on $[0, C_i h^i]$, where $C_i$ are defined by (\ref{Constante_01_2})-\ref{Constante_02_2}). \\[0.1cm]
Then, the probability of the event $\D\left\{ X^{(m)}(h) \leq X^{(k)}(h)\right\}$ is given by:
\vspace{-0.2cm}
\begin{equation}\label{Nonlinear_Prob}
\D Prob\left\{ X^{(m)}(h) \leq X^{(k)}(h)\right\} = \left |
\begin{array}{ll}
\D \hs 1 - \frac{1}{2}\!\left(\!\frac{\!\!h}{h^*}\!\right)^{\!\!m-k} & \mbox{ if } \hs 0 < h \leq h^*, \\[0.5cm]
\D \hs \frac{1}{2}\!\left(\!\frac{h^*}{\!\!h}\!\right)^{\!\!m-k} & \mbox{ if } \hs h \geq h^*.
\end{array}
\right.
\end{equation}
\end{theorem}
\begin{prooff}
$\frac{}{}$ \sa
$\blacktriangleright$ Let us first consider a fixed value of $h$ such that $h < h^*$.\sa
In this case, $f_m(h) < f_k(h)$, or in other words, $0 < C_m h^m < C_k h^k$ and due to Bramble-Hilbert lemma (see Figure \ref{Puissance}), one must deal with the following inequalities:
\begin{equation}\label{X_Y}
X^{(k)}(h) \leq C_k h^k \mbox{ and } X^{(m)}(h) \leq C_m h^m < C_k h^k.
\end{equation}
Then, to compute the probability such that $X^{(m)}(h) \leq X^{(k)}(h)$, we consider inequalities (\ref{X_Y}) in the plane $(0; X^{(m)}(h), X^{(k)}(h))$, (see Figure \ref{Trapeze}) in which the two random variables belong to the rectangle $R_t$ defined on $[0, C_m h^m] \times [0, C_k h^k]$. \sa
Our purpose is to characterize the points in $R_t$ that satisfy $X^{(m)}(h) \leq X^{(k)}(h)$. Obviously, it only concerns the points which are above the bisector $X^{(k)}(h)=X^{(m)}(h)$, namely the points which belong to the trapezium $T_u$ (see figure \ref{Trapeze}) whose surface is given by:
%
\begin{figure}[h]
  \centering
  \includegraphics[width=10cm]{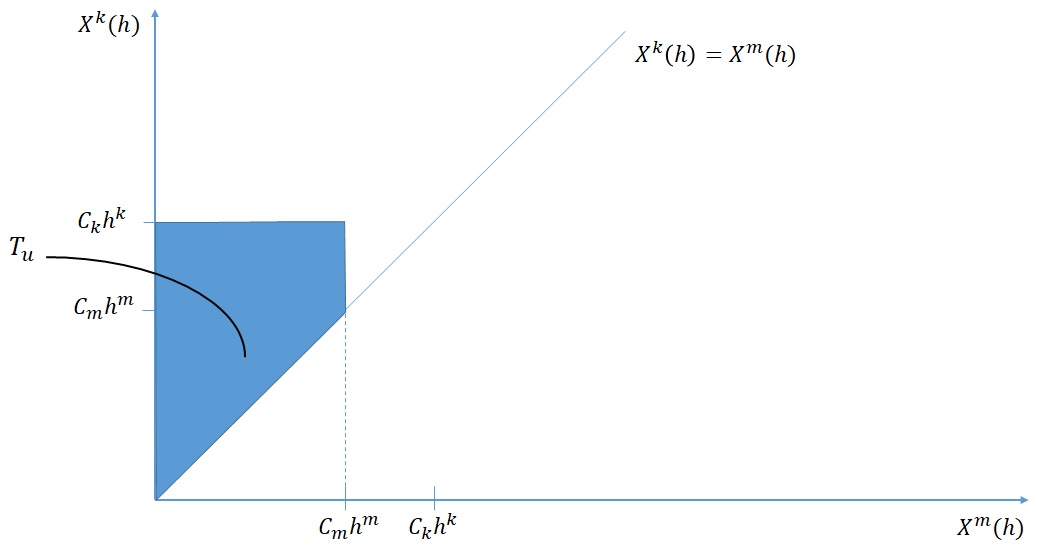}
  \caption{Area corresponding to $X^{(m)}(h) \leq X^k(h)$.} \label{Trapeze}
\end{figure}
\vspace{-0.2cm}
\begin{equation}
\D S(T_u) = C_m h^m (C_k h^k - C_m h^m) + \frac{C_m^{\,2}h^{\,2m}}{2},
\end{equation}
while the total surface of the rectangle $R_t$ is equal to $C_m C_k h^{m+k}$. \sa
As we assume that the two random variables $X^{(k)}(h)$ and $X^{(m)}(h)$ are independent and uniformly distributed, the probability $Prob\left\{ X^{(m)}(h) \leq X^{(k)}(h)\right\}$ corresponds to the ratio between the two surfaces of $T_u$ and $R_t$ and we have:
\begin{eqnarray}
\D Prob\left\{ X^{(m)}(h) \leq X^{(k)}(h)\right\} & = & \frac{S(T_u)}{S(R_t)} = \frac{C_m h^m (C_k h^k - C_m h^m) + C_m^{\,2}h^{\,2m}/2}{C_m C_k h^{m+k}}\,, \nonumber \\
& = & 1 - \frac{1}{2}\frac{C_m}{C_k}\,h^{m-k}.\label{Proba_1}
\end{eqnarray}
Using the definition (\ref{h*}) of $h^*\!$, we get:
\vspace{-0.4cm}
\begin{equation}
\D \forall h<h^*: Prob\left\{ X^{(m)}(h) \leq X^{(k)}(h)\right\} = 1 - \frac{1}{2}\frac{C_m}{C_k}\,h^{m-k} = 1 - \frac{1}{2}\!\left(\!\frac{\!\!h}{h^*}\!\right)^{\!\!m-k}\hspace{-0.4cm}.
\end{equation}
%
$\blacktriangleright$ Let us consider now the second case where $h>h^*$. \\ The curve $f_m(h)=C_m h^m$ is above the curve $f_k(h)=C_k h^k$ and by the same arguments we used above, one must deal with the following inequalities:
\begin{equation}\label{X_Y_1}
X^{(m)}(h) \leq C_m h^m \mbox{ and } X^{(k)}(h) \leq C_k h^k < C_m h^m.
\end{equation}
Then, if we change the role between $k$ and $m$, we can directly write :
\begin{eqnarray}
\D Prob\left\{X^{(k)}(h) \leq X^{(m)}(h)\right\} & = & \frac{C_k h^k (C_m h^m - C_k h^k) + C_k^{\,2}h^{\,2k}/2}{C_m C_k h^{m+k}} \nonumber \\
& = & 1 - \frac{1}{2}\,\frac{C_k}{C_m}h^{k-m}.\label{Proba_2}
\end{eqnarray}
Hence, the probability of the complementary event $X^{(m)}(h) \leq X^{(k)}(h)$ which interests us is given by:
\begin{eqnarray}
\D Prob\left\{X^{(m)}(h) \leq X^{(k)}(h)\right\} & = & 1- \D Prob\left\{X^k(h) \leq X^{(m)}(h)\right\} \nonumber \\
& = & \frac{1}{2}\frac{C_k}{C_k}.\frac{1}{h^{m-k}}  = \frac{1}{2}\!\left(\!\frac{h^*}{\!\!h}\!\right)^{\!\!m-k},
\end{eqnarray}
where we used the definition (\ref{h*}) of $h^*$.
\end{prooff}
The global shapes of the two probabilistic distributions (\ref{Heaviside_Prob}) and (\ref{Nonlinear_Prob}) are plotted in Figure \ref{Sigmoid} and particular features of (\ref{Nonlinear_Prob}) are described in the next section.
\section{Properties of the sigmoid probability distribution}\label{probabilistic_law_properties}
\noindent We give now the main properties of the sigmoid probability distribution given by (\ref{Nonlinear_Prob}). To this end, we will denote by $\mathcal{P}(h)$ the probability defined by:
\begin{equation}\label{P(h)}
\D \mathcal{P}(h) \equiv Prob\left\{X^{(m)}(h) \leq X^{(k)}(h)\right\}. \sa
\end{equation}
\begin{itemize}
\item The first feature we observe concerns the global shape of P(h) together with (\ref{Nonlinear_Prob}), drawn in Figure \ref{Sigmoid} for $m-k\neq 1$, which looks like a kind of sigmoid roughly approximated by a stepwise function given by (\ref{Heaviside_Prob}) from lemma \ref{Prob_General} of subsection \ref{two_steps}. \sa
In this way, we achieve our objective to relax the dependency assumption between the events $A$ and $B$. As a consequence, non linearity appears in the relation described by (\ref{Nonlinear_Prob}) between the probability of the event "\emph{$P_m$ finite element is more accurate than $P_k$ finite element}" and the mesh size $h$. \vspace{0.2cm}
\item \underline{Behavior of $\mathcal{P}(h)$ in the neighborhood of $0^+$}: \sa
Directly, we get:
\begin{equation}\label{limit_p(h)_en_0}
\D \lim_{h\rightarrow 0^+}\mathcal{P}(h) = \lim_{h\rightarrow 0^+}Prob\left\{ X^{(m)}(h) \leq X^{(k)}(h)\right\} = 1,
\end{equation}
which corresponds to the classical understanding of the error estimate (\ref{estimation_error}) which derives from Bramble-Hilbert lemma, namely asymptotically when the maximum of the mesh size $h$ goes to zero.\sa
Indeed, in these cases $h$ "\emph{is sufficiently small}", and despite the unknown values of the constants $C_k$ and $C_m$ which appear in (\ref{Constante_01_2}) and (\ref{Constante_02_2}), one concludes as expected that the finite element $P_m$ is more accurate than the finite element $P_m$, if $k<m$.\sa
But, the question is to determine what does it mean when $h$ "\emph{is sufficiently small}". We will partially discuss about this in the next point regarding the behavior of $\mathcal{P}(h)$ at the neighborhood of $h^*$ given by (\ref{h*}).\sa
From a probabilistic point of view the result (\ref{limit_p(h)_en_0}) is also intuitive because, when $h$ goes to $0^+$, the quantity $C_m h^m$ goes to 0 faster than $C_k h^k, (k<m)$. Depicting the relative position of $X^{(m)}(h)$ and $X^{(k)}(h)$ in a one dimensional way, (see Figure \ref{Axe}), it is clear that the probability of the event $\left\{X^{(m)}\!\leq X^{(k)}\right\}$ goes to 1 when $h$ goes to zero, as $X^{(m)}\! \leq C_m h^m$ due to Bramble-Hilbert lemma. \sa
\begin{figure}[h]
  \centering
  \includegraphics[width=8cm]{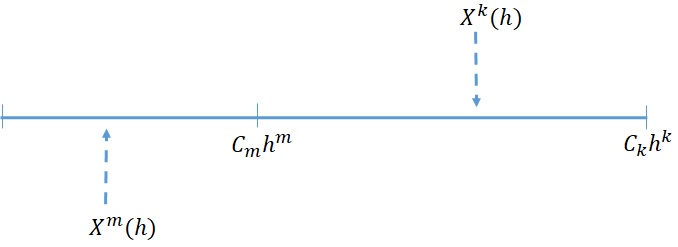}
  \caption{Relative one dimensional position between $X^{(m)}$ and $X^{(k)}$ ($h<h^*$).} \label{Axe}
\end{figure}
However, the interest of any probability distribution is to get additional information concerning the relative accuracy between two given finite elements, not only when $h$ goes to zero, as we will see further. Here, we just mentioned that we find again the well known conclusion to compare two finite elements when the mesh size is arbitrarily small.\sa
Indeed, finite element $P_m$ is not only {\em asymptotically} more accurate than $P_k$ as $k~<~m$. Indeed, for all $h \leq h^*$, the probability for $P_m$ to be more accurate than $P_k$ is between 0.5 to 1. It means that $P_m$ is {\em more likely} accurate than $P_k$ for all of these values of $h$. We also notice that we have not anymore the event "$P_m$ is more accurate than $P_k$" as an almost sure event as we got in subsection \ref{two_steps} with the law (\ref{Heaviside_Prob}). This is because we dropped the hypothesis of dependency between the events $A$ and $B$ which leads to a more general and realistic probabilistic distribution.\vspace{0.2cm}
\item \underline{Behavior of $\mathcal{P}(h)$ in the neighborhood of $h^*$}: \sa
The probabilistic stepwise law (\ref{Heaviside_Prob}) did not described the case $h$ equals $h^*$. However, here, the sigmoid probability distribution (\ref{Nonlinear_Prob}) can be extended by continuity to $h=h^*$ as we simply have:
\begin{equation}\label{limit_p(h)_en_h*_0}
\D \lim_{h\rightarrow h^{*-}}\mathcal{P}(h) = \hspace{-0.05cm}\lim_{h\rightarrow h^{*+}}\mathcal{P}(h) = \frac{1}{2},
\end{equation}
and then, we extend $\mathcal{P}(h)$ by continuity at $h^*$ by setting:
\begin{equation}
\mathcal{P}(h^*) = Prob\left\{X^{(m)}(h^*) \leq X^k(h^*)\right\} \equiv \frac{1}{2}.
\end{equation}
This feature illustrates that when $h=h^*$, $C_k h^{*k} = C_m h^{*m}$, and the two norms $\|u^{(k)}_h-u\|_{1,\Omega}$ and $\|u^{(m)}_h-u\|_{1,\Omega}$, which measures each approximation error of the two corresponding Lagrange finite elements, are somewhere below the two curves (see Figure \ref{Puissance}), or in other words, somewhere in the same interval as we here: $[0,C_kh^{*k}]=[0,C_mh^{*m}]$. Then, the probability to get $\D\left\{X^{(m)}(h^*) \leq X^k(h^*)\right\}$ is equal to 0.5.\sa
This new behavior claims that when $h$ approaches the critical value $h^*$ the event "\emph{$P_m$ finite element is more accurate than $P_k$ finite element}" is equally likely to occur or not to occur. As a consequence the accuracy between the two finite element $P_k$  and $P_m$ is equivalent. \sa
It is clearly a new theoretical information because, as we mentioned above, the values of the two constants $C_k$ and $C_m$ are totally unknown. Indeed, we already suspected and pointed out by data mining techniques, (see for example \cite{AsCh11}, \cite{AsCh13} and \cite{AsCh14}), that this situation would occur. Here, we complete this suspicion by a theoretical probabilistic framework.\vspace{0.2cm}
\item Despite the usual point of view which claims that $P_m$ finite element are more accurate than $P_k$ ones, we get here that $P_k$ finite element is \emph{more likely} accurate than $P_m$ when $h>h^*$. This new point of view allows us to recommend that for specific situations, like for adaptive refinement meshes for example, $P_k$ finite element would be locally more appropriated as long as one will be able to detect the case $h>h^*$.
\end{itemize}
\vspace{-0.4cm}
\section{Conclusions}
\noindent In this paper we present a new way to investigate the relative accuracy between two finite elements. Indeed, leaving the classical asymptotic point of view usually considered to compare the speed of convergence for different approximation errors, we got new insights for understanding error estimates. The way we thought the error estimates is not restricted to the finite element method but can be extended to other approximation methods. Indeed, the underlying idea is that, given a class of numerical schemes and their corresponding error estimates, one is able to rank them, not only in terms of asymptotic speed of convergence as usual, but also by evaluating the \emph{almost surely} more accurate one. \sa
For example, considering numerical schemes to approximate solution to ordinary differential equations, one would be able to argue, why (or why not!) RK4 scheme would be implemented rather than another simplest one. \sa
\textbf{\underline{Homages}:} The authors want to warmly dedicate this research to pay homage to the memory of Professors Andr\'e Avez and G\'erard Tronel who largely promote the passion of research and teaching in mathematics.

\end{document}